\documentclass[10pt]{amsart}

\newtheorem{thm}{Theorem}
\newtheorem{lemma}{Lemma}

\newtheorem{cor}{Corollary}
\newtheorem{prop}{Proposition}

\usepackage{amsmath}
\usepackage[psamsfonts]{amssymb}
\usepackage[mathscr]{euscript}

\usepackage[mathscr]{euscript}

\newcommand{\cal}{\mathcal}

\title[\hfill]{Curvature properties of twistor  spaces}

\author{Johann Davidov and Oleg Mushkarov}
\thanks{The  authors are partially supported by  the National Science
Fund, Ministry of Education and Science of Bulgaria under contract
DN 12/2 }

\dedicatory{Dedicated to Professor Armen Sergeev on the occasion
of his 70th birthday}

\address{Johann Davidov\\Institute of Mathematics and Informatics \\
Bulgarian Academy of Sciences\\ Acad. G.Bonchev Str. Bl.8\\ 1113
Sofia\\ Bulgaria} {\email{jtd@math.bas.bg}

\address{Oleg Mushkarov \\Institute of Mathematics and Informatics \\
Bulgarian Academy of Sciences\\ Acad. G.Bonchev Str. Bl.8\\ 1113
Sofia\\ Bulgaria\\ \newline \centerline{and} \newline Sauth-West
University, 2700 Blagoevgrad, Bulgaria} \email{muskarov@math.bas.bg}

\begin{document}

\begin{abstract}  In this paper
we review some results on the Riemannian and almost Hermitian
geometry of twistor spaces of oriented Riemannian $4$-manifolds with
emphasis on their curvature properties.

\vspace{0,1cm} \noindent 2010 {\it Mathematics Subject
Classification}. Primary 53C28, Secondary 53C56.

\vspace{0,1cm} \noindent {\it Key words: twistor spaces, almost
Hermitian structures, Riemannian and Hermitian curvatures of twistor
spaces}.
\end{abstract}

\thispagestyle{empty}

\maketitle \vspace{0.5cm}

\section{Introduction}

The twistor theory has its origin in Mathematical Physics. Inspired
by the Penrose program \cite{Pe1,Pe2}, Atiyah, Hitchin and Singer
\cite{AHS} developed this theory on oriented Riemannian
$4$-manifolds. They defined the  twistor space of such a manifold
$M$ as the $2$-sphere bundle ${\mathcal Z}$ on $M$ whose fiber at
any point $p\in M$ consists of all complex structures on the tangent
space $T_pM$ compatible with the metric and the opposite orientation
of $M$. The $6$-manifold ${\mathcal Z}$ admits two natural almost
complex structures ${\cal J}_1$ and ${\cal J}_2$ introduced,
respectively, by Atiyah, Hitchin, and Singer \cite{AHS} and by Eells
and Salamon \cite{ES85}. The almost complex structure ${\cal J}_2$
is never integrable but it plays an important role in the theory of
harmonic maps. The almost complex structure  ${\cal J}_1$ is
integrable if and only if the base manifold $M$ is self-dual, i.e.
its Weyl conformal tensor ${\mathcal W}$ satisfies the equation
$\ast {\mathcal W}={\mathcal W}$, where $\ast$ is the Hodge
star-operator \cite{Besse}. So, in the case when $M$ is self-dual,
$\mathcal Z$ is a complex $3$-manifold and each fibre $\mathcal
Z_p=\pi^{-1}(p), p\in M$,  is a complex submanifold of $\mathcal Z$
biholomorphic to $\mathbb{C}\mathbb{P}^1$. The normal bundle of
$\mathcal Z_p$ is biholomorphically equivalent to $H\oplus H$, where
$H$ is the positive Hopf bundle on $\mathbb{C}\mathbb{P}^1$. The
antipodal map $j\rightarrow -j$ on each fibre induces an
anti-biholomorphic involution of $\mathcal Z$ without fixed points.

The above construction depends only on the conformal class of the
given metric on $M$ and, conversely, the complex structure of the
twistor space determines the self-dual conformal structure on $M$.
More precisely, let $\mathcal Z$ be a complex $3$-manifold with the
following properties: (1) $\mathcal Z$ is fibered by projective
lines whose normal bundle is isomorphic to $H\oplus H$ ; (2)
$\mathcal Z$ possesses a free anti-holomorphic involution which
transforms each fibre to itself. Then $\mathcal Z$ is the twistor
space of a self-dual manifold $M$ (\cite {AHS}, see also
\cite{Besse}). The described correspondence between self-dual
manifolds and twistor spaces is called the Penrose correspondence
and it has been used for years to study the conformal geometry of
four-manifolds by means of complex geometry techniques.

The twistor space $\mathcal Z$ admits a $1$-parameter family of
Riemannian metrics which are compatible with the almost complex
structures ${\cal J}_1$ and ${\cal J}_2$  and the natural projection
$\pi : {\mathcal Z}\to M$ is a Riemannian submersion \cite{Besse}.
These natural almost Hermitian structures are very interesting
geometric objects in their own right whose differential geometric
properties have been studied by many authors, to cite here just a
few \cite{AGI98, De, FD09, FG85, FK82, G, Hit81, Hit95, Pon}.

Motivated by  open questions in  non-K\"ahler geometry in a series
of papers \cite{M87, DM, DM90, DGM01, DGM02, DMG95, DM18, ADM,
ADM14}, the authors studied the twistor spaces of oriented
Riemannian $4$-manifolds as a source of examples of almost Hermitian
$6$-manifolds with interesting geometric properties. In the present
paper we review some of these results with emphasis on the curvature
properties of twistor spaces.

\medskip

{\Large{Table of contents}}

\smallskip

2. Preliminaries

3. Riemannian curvature of twistor spaces

   \ \ \ \ 3.1. Riemannian sectional curvature

   \ \ \ \ 3.2. Ricci curvature

4. Twistor spaces with Hermitian Ricci tensor

5. K\"ahler curvature identities on twistor spaces

6. $\ast$-Einstein twistor spaces

7. Curvature properties of the Chern connection on twistor spaces

8. Holomorphic curvatures of twistor spaces

   \ \ \ \ 8.1. Holomorphic bisectional curvature

   \ \ \ \ 8.2. Orthogonal bisectional curvatures

   \ \ \ \ 8.3. Hermitian bisectional curvature

\medskip

\noindent {\bf Acknowledgements}.  The authors would like to thank
the editors of this volume for the invitation to submit a paper in
honour of  Professor Armen Sergeev, our dear friend and colleague.

\section{Preliminaries}

Let $(M,g)$ be an oriented (connected) Riemannian manifold of
dimension four. The metric $g$ induces a metric on the bundle of
two-vectors $\pi:\Lambda^2TM\to M$ by the formula
$$
g(v_1\wedge v_2,v_3\wedge v_4)=\frac{1}{2}det[g(v_i,v_j)],
$$
the  factor $1/2$ being chosen in accordance  with \cite{DM}. Let
$\ast:\Lambda^kTM\to \Lambda^{4-k}M$, $k=0,...,4$, be the Hodge star
operator. For $k=2$, it is an involution of $\Lambda^2TM$, thus we
have the orthogonal decomposition
$$
\Lambda^2TM=\Lambda^2_{-}TM\oplus\Lambda^2_{+}TM,
$$
where $\Lambda^2_{\pm}TM$ are the subbundles of $\Lambda^2TM$
corresponding to the eigenvalues $\pm 1$ of the operator $\ast$.
Given a local oriented orthonormal frame $(E_1,E_2,E_3,E_4)$ of $TM$
we set
\begin{equation}\label{s-basis}
s_1^{\pm}=E_1\wedge E_2\pm E_3\wedge E_4, \quad s_2^{\pm}=E_1\wedge
E_3\pm E_4\wedge E_2, \quad s_3^{\pm}=E_1\wedge E_4\pm E_2\wedge
E_3.
\end{equation}
Then $(s_1^{\pm},s_2^{\pm},s_3^{\pm})$ is a local orthonormal frame
of $\Lambda^2_{\pm}TM$. This frame defines an orientation on
$\Lambda^2_{\pm}TM$ which does not depend on the choice of the frame
$(E_1,E_2,E_3,E_4)$ (see, for example, \cite{D17}). We call this
orientation "canonical".

For every $a\in\Lambda ^2TM$, define a skew-symmetric endomorphism
of $T_{\pi(a)}M$ by

\begin{equation}\label{cs}
g(K_{a}X,Y)=2g(a, X\wedge Y), \quad X,Y\in T_{\pi(a)}M.
\end{equation}
Denoting by $G$ the standard metric $-\frac{1}{2}Trace\,(PQ)$ on the
space of skew-symmetric endomorphisms, we have $G(K_a,K_b)=2g(a,b)$
for $a,b\in \Lambda ^2TM$. The assignment $a\to K_a$ is the standard
isomorphism of the bundle $\Lambda ^2TM$ with the bundle $A(TM)$ of
$g$-skew-symmetric endomorphism of $TM$. Let $\times$ be the usual
vector cross product on the oriented $3$-dimensional vector space
$\Lambda^2_{\pm}T_pM$, $p\in M$, endowed with the metric $g$. If
$a,b\in\Lambda^2_{\pm}T_pM$, the isomorphism $\Lambda^2TM\cong
A(TM)$ sends $a\times b$ to $\pm\frac{1}{2}[K_a,K_b]$. In the case
when $a\in\Lambda^2_{+}T_pM$, $b\in\Lambda^2_{-}T_pM$, the
endomorphisms $K_a$ and $K_b$ of $T_pM$ commute. For
$a,b\in\Lambda_{\pm}T_pM$, we have
$$
K_{a}\circ K_{b}=-g(a,b)Id\pm K_{a\times b}.
$$
In particular, $K_a$ and $K_b$, $a,b\in\Lambda_{\pm}T_pM$,
anti-commute if and only if $a$ and $b$ are orthogonal.

If $\sigma\in\Lambda^2TM$ is of unit length, then $K_{\sigma}$ is a
complex structure on the vector space $T_{\pi(\sigma)}M$ compatible
with the metric $g$, i.e., $g$-orthogonal. Conversely, the
$2$-vector $\sigma$ dual to one half of the K\"ahler $2$-form of
such a complex structure is a unit vector in $\Lambda^2TM$.
Therefore the unit sphere bundle ${\mathcal Z}$ of $\Lambda ^2TM$
parametrizes the complex structures on the tangent spaces of $M$
compatible with the metric $g$. This bundle is called the twistor
space of the Riemannian manifold $(M,g)$. Since $M$ is oriented, the
manifold ${\mathcal Z}$ has two connected components ${\mathcal
Z}_{\pm}$ called the positive and  negative twistor spaces of
$(M,g)$. These are the unit sphere subbundles of
$\Lambda^2_{\pm}TM$. The bundle $\pi: {\mathcal Z}_{\pm}\to M$
parametrizes  the complex structures on the tangent spaces of $M$
compatible with the metric and  the  $\pm$ orientation via the
correspondence ${\mathcal Z}_{\pm}\ni\sigma\to K_{\sigma}$.

The connection $\nabla$ on $\Lambda^2TM$ induced by the Levi-Civita
connection of $M$ preserves the bundles $\Lambda^2_{\pm}TM$, so it
induces a metric connection on each of them denoted again by
$\nabla$. The horizontal distribution of $\Lambda^2_{\pm}TM$ with
respect to $\nabla$ is tangent to the twistor space ${\mathcal
Z}_{\pm}$. Thus, we have the decomposition $T{\mathcal
Z}_{\pm}={\mathcal H}\oplus {\mathcal V}$ of the tangent bundle of
${\mathcal Z}_{\pm}$ into horizontal and vertical components. The
vertical space ${\mathcal V}_{\sigma}=\{V\in T_{\sigma}{\mathcal
Z}_{\pm}:~ \pi_{\ast}V=0\}$ at a point $\sigma\in{\mathcal Z}_{\pm}$
is the tangent space to the fibre of ${\mathcal Z}_{\pm}$ through
$\sigma$. If we consider $T_{\sigma}{\mathcal Z}_{\pm}$ as a
subspace of $T_{\sigma}(\Lambda^2_{\pm}TM)$, then the space
${\mathcal V}_{\sigma}$ is the orthogonal complement of $\sigma$ in
$\Lambda^2_{\pm}T_{\pi(\sigma)}M$.

The $6$-manifold ${\mathcal Z}_{\pm}$ admits two almost complex
structures ${\cal J}_1$ and ${\cal J}_2$ introduced, respectively,
by Atiyah, Hitchin, and Singer \cite{AHS} and by Eells and Salamon
\cite{ES85}. For $\sigma\in{\mathcal Z}_{\pm}$, the horizontal space
${\mathcal H}_{\sigma}$ is isomorphic to the tangent space
$T_{\pi(\sigma)}M$ via the differential $\pi_{\ast\,\sigma}$, and
the structures ${\mathcal J}_1$ and ${\mathcal J}_2$  on the space
${\mathcal H}_{\sigma}$ are both defined as the lift to ${\mathcal
H}_{\sigma}$ of the complex structure $K_{\sigma}$ on
$T_{\pi(\sigma)}M$. The vertical space ${\mathcal V}_{\sigma}$ is
tangent to the unit sphere in the $3$-dimensional vector space
$(\Lambda^2_{\pm}T_{\pi(\sigma)}M,g)$, and we denote by ${\mathcal
J}_{\sigma}$ the standard complex structure of the unit sphere
restricted to ${\mathcal V}_{\sigma}$. It is given by
$$
{\mathcal J}_{\sigma}V=\pm(\sigma\times V),\quad V\in{\mathcal
V}_{\sigma},
$$
 On a vertical space ${\mathcal V}_{\sigma}$,
${\mathcal J}_1$ is defined to be the complex structure ${\mathcal
J}_{\sigma}$ of the fibre through $\sigma$, while ${\mathcal J}_2$
is defined as the conjugate complex structure, i.e., ${\mathcal
J}_{2}|{\mathcal V}_{\sigma}=-{\mathcal J}_{\sigma}$. Thus, for
$\sigma\in{\mathcal Z}_{\pm}$,
\begin{equation}\label{ACS}
\begin{array}{c}
{\mathcal J}_{n}|{\mathcal H}_{\sigma}=(\pi_{\ast}|{\mathcal
H}_{\sigma})^{-1}\circ K_{\sigma}\circ\pi_{\ast}|{\mathcal
H}_{\sigma}\\[8pt]

{\cal J}_{n}V=\pm(-1)^{n+1}(\sigma\times V )\;~~\mbox {for} ~~\;
V\in \mathcal{V_{\sigma}},\;\;\;\;n=1,2.
\end{array}
\end{equation}

It is a result of Eells and Salamon \cite{ES85} that the almost
complex structure ${\mathcal J}_2$ is never integrable, so it does
not come from a complex structure. Nevertheless, ${\mathcal J}_2$ is
very useful for constructing harmonic maps. The integrability
condition for ${\mathcal J}_1$ has been found by Atiyah, Hitchin,
and Singer \cite{AHS}. To state their result, we first recall the
well-known curvature decomposition in dimension four. Note that for
the curvature tensor $R$ , we adopt the following definition:
$R(X,Y)=\nabla_{[X,Y]}-[\nabla_{X},\nabla_{Y}]$. The curvature
operator ${\mathcal R}$ corresponding to the curvature tensor is the
endomorphism of $\Lambda ^{2}TM$ defined by
$$
g({\mathcal R}(X\wedge Y),Z\wedge U)=g(R(X,Y)Z,U),\quad X,Y,Z,U\in
TM.
$$
Denote by $\rho$ the Ricci tensor of $(M,g)$ and by  $A:TM\to TM$
its Ricci operator, $g(A(X),Y)=\rho(X,Y)$. Then the endomorphism
${\cal B}:\Lambda^2TM\to \Lambda^2TM$ corresponding to the
traceless Ricci tensor is given by
\begin{equation}\label{B}
{\mathcal B}(X\wedge Y)=A(X)\wedge Y+X\wedge
A(Y)-\frac{s}{2}X\wedge Y,
\end{equation}
where $s$ is the scalar curvature. Note that ${\mathcal B}$ sends
$\Lambda^2_{\pm}TM$ into $\Lambda^2_{\mp}TM$. Let ${\mathcal W}:
\Lambda^2TM\to \Lambda^2TM$ be the endomorphism corresponding to
the Weyl conformal tensor. Denote the restriction of ${\mathcal
W}$ to $\Lambda^2_{\pm}TM$ by ${\mathcal W}_{\pm}$, so ${\mathcal
W}_{\pm}$ sends $\Lambda^2_{\pm}TM$ to $\Lambda^2_{\pm}TM$ and
vanishes on $\Lambda^2_{\mp}TM$.

It is well known that the curvature operator decomposes as
(\cite{ST69}, \cite[Chapter 1 H]{Besse})
\begin{equation}\label{dec}
{\mathcal R}=\frac{s}{6}Id+{\mathcal B}+{\mathcal W}_{+}+{\mathcal
W}_{-}
\end{equation}
Note that this differs from \cite{Besse} by  a factor of $1/2$
because of the factor $1/2$ in our definition of the induced metric
on $\Lambda^2TM$. The Riemannian manifold $(M,g)$ is Einstein
exactly when ${\mathcal B}=0$. It is called self-dual
(anti-self-dual) if ${\mathcal W}_{-}=0$ (${\mathcal W}_{+}=0$). The
self-duality (anti-self-duality) condition is invariant under
conformal changes of the metric since the Weyl tensor is so. Note
also that changing the orientation of $M$ interchanges the roles of
$\Lambda^2_{-}TM$ and $\Lambda^2_{+}TM$ (respectively, of ${\mathcal
Z}_{-}$ and ${\mathcal Z}_{+}$), hence the roles of ${\mathcal
W}_{-}$ and ${\mathcal W}_{+}$.

The famous Atiyah-Hitchin-Singer  theorem \cite{AHS} states that
the almost complex structure ${\cal J}_1$ on ${\cal Z}_{-}$ (resp.
${\cal Z}_{+}$) is integrable if and only if $(M,g)$ is self-dual
(resp. anti-self-dual).

The twistor space ${\cal Z}_{\pm}$ of an oriented Riemannian
$4$-manifold $(M,g)$ admits a natural $1$-parameter family of
Riemannian metrics $h_t$ defined by
$$
h_{t}=\pi ^{*}g+tg^{v}
$$
where $t>0$ and $g^{v}$ is  the restriction of the metric of
$\Lambda^{2}TM$ on the vertical distribution ${\cal V}$. Then
$\pi:({\cal Z}_{\pm},h_t)\to (M,g)$ is a Riemannian submersion
with totally geodesic fibres, and the almost-complex structures
$\cal{J}_{1}$ and $\cal{J}_{2}$ are compatible with the metrics
$h_{t}$. The Gray-Hervella classes \cite{GH} of the almost
Hermitian structures $(h_t,\cal{J}_n), t>0, n=1,2$, have been
determined in \cite{M87}.

\section{Riemannian curvature of twistor spaces}

The O'Neill formulas \cite{O'N}, \cite[Ch.\,9 G]{Besse} can be used
to obtain coordinate-free formulas for  various curvatures of the
metric $h_t$ on the twistor space in terms of the curvature of its
base manifold $M$. This is done in \cite{DM} in the case when
$dim\,M=4$ and in \cite{D05} for the general twistor space of
partially complex structures ($f$-structures) on a Riemannian
manifold of any dimension $\geq 3$. We shall discuss here the most
interesting case of the negative twistor space of an oriented four
dimensional  Riemannian manifold. The reason to choose the negative
twistor space is connected with the Atiyah-Hitchin-Singer theorem
mentioned above. As smooth manifolds, the positive and the negative
twistor spaces of ${\mathbb C}{\mathbb P}^2$ coincide with the
complex flag manifold $F_{1,2}$. The Atiyah-Hitchin-Singer almost
complex structure on the negative twistor space of ${\mathbb
C}{\mathbb P}^2$ is integrable and coincides with the standard
complex structure of $F_{1,2}$, while it is not integrable  on the
positive twistor space.

\smallskip

In what follows, $(M,g)$ will denote an oriented Riemannian
manifold of dimension four, and ${\mathcal Z}$ will stand for its
negative twistor space ${\mathcal Z}_{-}$.

\smallskip

\subsection{Riemannian sectional curvature}

Let $(M,g)$ be an oriented Riemannian $4$-manifold with Levi-Civita
connection  $\nabla$ and Riemannian curvature tensor $R$. For any
$t>0$ denote by $R_t$ the Riemannian curvature tensor of the metric
$h_t$ on the twistor space ${\mathcal Z}$ of $(M,g)$. Applying the
O'Neill formulas \cite{O'N} for the Riemannian submersion
$\pi:({\mathcal Z},h_t)\to (M,g)$, one can compute the Riemannina
sectional curvature of $({\mathcal Z},h_t)$.

\begin{prop}\label{Sec}\rm (\cite{DM}) Let $E,F\in T_{\sigma}{\mathcal Z}$,
$X=\pi_{\ast}E, Y=\pi_{\ast}F,V={\mathcal V}E$ and $ W={\mathcal
V}F$. Then
$$
\begin{array}{c}\label{Sec}
h_t(R_t(E\wedge F)E,F)=g(R(X\wedge Y)X,Y)-tg((\nabla_X\mathcal{R})(X\wedge Y),\sigma\times W)\\[6pt]
+tg((\nabla_Y\mathcal{R})(X\wedge Y),\sigma\times
V)-3tg(\mathcal{R}(\sigma),X\wedge
Y)g(\sigma\times V,W)\\[6pt]
-t^2g(R(\sigma\times V)X,R(\sigma\times
W)Y)+\displaystyle{\frac{t^2}{4}}\|R(\sigma\times W)X+R(\sigma\times
V)Y\|^2\\[6pt]
-\displaystyle{\frac{3t}{4}}\|R(X\wedge
Y)\sigma\|^2+t(\|V\|^2\|W\|^2-g(V,W)^2).
\end{array}
$$
\end{prop}
In the case when the base manifold $(M,g)$ is self-dual and Einstein the above
formula takes an apparently simple form.

\begin{cor}\label{SecS} Let $(M,g)$ be a  self-dual Einstein
manifold with scalar curvature $s$. Then
$$ h_t(R_t(E\wedge F)E,F)=g(R(X\wedge Y)X,Y)-\displaystyle{\frac{ts}{2}}g(\sigma,X\wedge Y)g(\sigma \times
V,W)$$
$$-\displaystyle{\frac{1}{2}(\frac{ts}{12}})^2g(X,Y)g(V,W)+3\displaystyle{(\frac{ts}{12}})^2g(X\wedge Y,V\wedge W)$$
$$+\displaystyle{(\frac{ts}{24}})^2(\|X\|^2\|W\|^2+\|Y\|^2\|V\|^2)$$
$$-6t\displaystyle{(\frac{s}{24}})^2(\|X\wedge Y\|^2-2g(\sigma,X\wedge Y)^2)$$
$$+t(\|V\|^2\|W\|^2-g(V,W)^2).$$
\end{cor}

\smallskip

\subsection {Ricci curvature}
\label{Ricci}

\bigskip

The study of the Ricci curvature of a twistor space is based on
the following explicit formula for the Ricci tensor which is a
consequence of Proposition \ref{Sec}.


\begin{prop}\label {Ricci}\rm (\cite{DM})
Let $\rho_{t}$ be the Ricci tensor of the twistor space
$({\mathcal Z},h_t)$. If $E\in T_{\sigma}{\mathcal Z}$,
$X=\pi_{\ast}E$, and $V={\mathcal V}E$, then
$$
\begin{array}{c}\label{Ricci}
\rho_{t}(E,E)=\rho(X,X)+tg(\delta{\mathcal R}(X),\sigma\times
V)+\displaystyle{\frac{t^2}{4}}||{\mathcal R}(\sigma\times
V)||^2\\[6pt]
+\displaystyle{\frac{t}{2}}||\imath_{X}\circ {\mathcal
R}(\sigma)||^2-\displaystyle{\frac{t}{2}}||\imath_{X}\circ {\mathcal
R}_-||^2 +||V||^2,
\end{array}
$$
where $\rho$ is the Ricci tensor of $(M,g)$, $\delta{\mathcal R}$
is the co-differential of ${\mathcal R}$, ${\mathcal R}_-$ is the
restriction of ${\mathcal R}$ on $\Lambda^2_{-}TM$, and
$\imath_{X}: \Lambda^2TM\to TM$ is the interior product.
\end{prop}

Taking the trace of $\rho_{t}$, we obtain the following formula for
the scalar curvature $s_t$ of the twistor space $({\mathcal
Z},h_t)$.

\begin{cor}\label {scalar}\rm (\cite{DM})
Let $s$ be the scalar curvature of $(M,g)$. Then
$$s_t(\sigma)=s(\pi(\sigma))+\frac{t}{4}(||{\mathcal
R}(\sigma)||^2-||{\mathcal R}_-||_{\pi(\sigma)}^2)+ \frac{2}{t}.
$$
\end{cor}

In the case when the base manifold of a twistor space is Einstein
and self-dual, these  formulas  can significantly be simplified as
follows.
\begin{cor}\label{R-ES}\rm (\cite{DM})
If $(M,g)$ is Einstein and self-dual, the Ricci tensor $\rho_t$ of
$({\mathcal Z},h_t)$ and its scalar curvature $s_{t}$ are given by
$$
\begin{array}{c}
\rho_{t}(E,E)=\displaystyle{[\frac{s}{4}-t(\frac{s}{12})^2]||X||^2+[1+(\frac{ts}{12})^2]||V||^2},\quad
E=X^h_{\sigma}+V,\\[10pt]
s_{t}=\displaystyle{\frac{2}{t}+s-\frac{t}{72}s^2}.
\end{array}
$$
\end{cor}

As an application of  Proposition~\ref{Ricci}, one can prove the
following result of T. Friedrich and R. Grunewald \cite{FG85} about
the Einstein condition on $({\mathcal Z},h_t)$.

\begin{thm}\label{Einstein}\rm (\cite{FG85,DM})
The Riemannian manifold $({\mathcal Z},h_t)$ is Einstein if and
only if $(M,g)$ is a self-dual Einstein manifold with scalar
curvature $s=6/t$ or $s=12/t$.
\end{thm}

The next  property of the Ricci tensor $\rho_t$ is an easy
consequence of Corollary \ref{R-ES}.
\begin{prop}\label{Ein-like}\rm (\cite{DMG92})
If $(M,g)$ is Einstein and self-dual, the covariant derivative of
the Ricci tensor $\rho_t$ of $({\mathcal Z},h_t)$ satisfies the
identity
\begin{equation}\label{E-l}
(D_{E}\rho_{t})(E,E)=0,\quad E\in T{\mathcal Z},
\end{equation}
where $D$ is the Levi-Civita connection of $h_t$. Moreover,
$\rho_t$ is parallel if and only if $st=6$, $st=12$, or $s=0$.
\end{prop}

\noindent {\bf Remark}. Condition (\ref{E-l}) for the Ricci tensor
$\rho$ of a Riemannian manifold $(N,h)$ is known as the third Ledger
condition \cite{L54}, \cite[Sec. 6.8]{Will}. It is easy to see by
polarization that (\ref{E-l}) is equivalent to the identity
$$
(\nabla_{X}\rho)(Y,Z)+(\nabla_{Y}\rho)(Z,X)+(\nabla_{Z}\rho)(X,Y)=0,
\ \  X,Y,Z\in TN,
$$
where $\nabla$ is the Levi-Chivita connection of $(N,h)$. If this
condition is satisfied, the manifold is real-analytic \cite{S93},
and the scalar curvature is constant \cite[Proposition 2.3]{DAN60}.
Condition ({\ref{E-l}) appears in the study of the so-called D'Atri
spaces which are characterized by the property that the geodesic
symmetries preserve the volume up to sign \cite{DAN60}. It is  one
of the Einstein-like conditions introduced and studied by A. Gray
\cite{G78}, and  discussed also in  Besse's book \cite[Sec.
16G]{Besse} as an interesting generalization of the Einstein
condition. We refer the reader to
\cite{Besse,D05,G78,Jel,KNa87,PedTod} for examples of Riemannian
manifolds satisfying condition (\ref{E-l}).
Proposition~\ref{Ein-like}, which gives twistorial examples of such
manifolds, seems to be interesting in the case of negative scalar
curvature of $(M,g)$ since the complete classification of compact
Einstein self-dual manifolds with negative scalar curvature is not
available yet. It has been conjectured by A. Vitter \cite{V86} that
every such a manifold is a quotient of the unit ball in ${\mathbb
C}^2$ with the metric of negative constant sectional curvature or
the Bergman metric.

Corollary~\ref{R-ES} can be used to show that an isometry of the
twistor space preserves vertical, and hence horizontal, spaces.
This implies the following.

\begin{lemma}\label{homog} \rm (\cite{D05})
If $(M,g)$ is an Einstein and self-dual manifold with scalar
curvature $s$, then every (local) isometry of the twistor space
$\pi: ({\mathcal Z},h_t)\to (M,g)$ descends to a (local) isometry
of the metric $g$ provided $ts\neq 6 $ and $ts\neq 12$.
\end{lemma}

\noindent {\bf Remarks}. \rm (\cite{D05}) Suppose that the manifold
$(M,g)$ is Einstein and self-dual, and $ts=6$ or $ts=12$.

\medskip

\noindent {\bf 1.} Lemma~\ref{homog} does not hold  as there may
exist an isometry of the twistor space of $(M,g)$ which does not
descend to an isometry of $g$. For example, it is well-known that
the twistor space ${\mathcal Z}$ of the sphere $S^4$ considered with
its standard metric is the complex projective space ${\Bbb C}{\Bbb
P}^3$. To describe the twistor projection $\pi:{\Bbb C}{\Bbb P}^3\to
S^4$, it is convenient to identify $S^4$ with the quaternionic
projective space ${\Bbb H}{\Bbb P}^1$. Then $\pi$ is given in
homogeneous coordinates by $[z_1,z_2,z_3,z_4]\to [z_1+z_2{\bf
j},z_3+z_4{\bf j}]$. If $ts=12$, the metric $h_t$ is a multiple of
the Fubini-Study metric. The map $\Psi:{\Bbb C}{\Bbb P}^3\to {\Bbb
C}{\Bbb P}^3$ defined by $\Psi([z_1,z_2,z_3,z_4])=[\frac{1}{\sqrt
2}(z_1+z_2),\frac{1}{\sqrt 2}(z_1-z_2),z_3,z_4]$ is an isometry of
the Fubini-Study metric which does not preserve all fibres of the
twistor projection $\pi$.

\medskip

\noindent {\bf 2.} The scalar curvature $s$ of $M$ is positive and,
by a result of Hitchin \cite{Hit81} and of Friedrich and Kurke
\cite{FK82}, see also \cite[Theorem 13.30]{Besse}, $(M,g)$ is
isometric to the sphere $S^4$ or the complex projective space ${\Bbb
C}{\Bbb P}^2$ with their standard metrics. In particular, the metric
$g$ is homogeneous, hence all of the metrics $h_t$ on the twistor
space are also homogeneous.

  The latter remark and Lemma~\ref{homog} give the following result
which seems to be  "folklore".

\begin{prop}\label{isometry} \rm (\cite{D05})
Let $(M,g)$ be a complete  Einstein self-dual  manifold. The metric
$h_t$ (with arbitrary $t$) on  the twistor space ${\mathcal Z}$ is
(locally) homogeneous if and only if the metric $g$ on the base
manifold $M$ is (locally) homogeneous.
\end{prop}

Proposition~\ref{Ein-like} and Lemma~\ref{homog} imply the
following
\begin{prop} \rm ({\cite{D05})}
Let $M$ be an inhomogeneous Einstein self-dual  $4$-manifold with
non-zero scalar curvature $s$. Then, for any $t>0$ with $ts\neq 6$
and $ts\neq 12$, the twistor space $({\mathcal Z},h_t)$ is
non-homogeneous, has non-parallel Ricci tensor satisfying the third
Ledger condition (\ref{E-l}) and is not locally isometric to a
Riemannian product.

Moreover, if $M$ is locally non-homogeneous, then so is its
twistor space.
\end{prop}

\noindent {\bf Remark}. \rm (\cite{D05}) If the base manifold $M$ is
locally homogeneous, so is its twistor space. There are a lot of
examples of (non-compact) locally non-homogeneous,  self-dual,
Einstein manifolds with non-zero scalar curvature, to cite just a
few papers where such examples (complete or not) can be
found:\cite{AG02,CP02,Der81,Hit95,LB91}.

\smallskip

\section{Twistor spaces with Hermitian Ricci
tensor}\label{Hermitian Ricci}

It is well-known that on any symplectic manifold $N$ with symplectic
form $\Omega$ there exist a Riemannian metric $h$ and a
$h$-orthogonal almost complex structure $J$ such that $\Omega$ is
the K\"ahler $2$-form of the almost Hermitian manifold $(N,h,J)$,
i.e., $\Omega(X,Y)=h(JX,Y)$ for $X,Y \in TN$. Recall that an almost
Hermitian manifold with closed K\"ahler $2$-form is called almost
K\"ahler. A Riemannian metric $h$ on $N$ is said to be associated to
the symplectic form $\Omega$ if there exists a $h$-orthogonal almost
complex structure $J$ for which $\Omega(X,Y)=h(JX,Y)$. Note that
such an almost complex structure is unique. Assume that $N$ is
compact, and denote by ${\mathcal A}$ the set of all Riemannian
metrics on $N$ associated to $\Omega$. If $h\in {\mathcal A}$ and
$J$ is the corresponding almost complex structure, let $s$ and
$s^{\ast}$ be the scalar curvature of the metric $h$ and the
$\ast$-scalar curvature of the almost Hermitian structure $(h,J)$
(we recall the definition of $s^{\ast}$ in Section 6). Then we can
consider the integrals
$$
\int_{N}s\,vol_h\quad {\rm{and}}\quad \int_{N}(s-s^{\ast})\,vol_h
$$
as functionals on the set ${\mathcal A}$. D. Blair and S. Ianu\c{s}
\cite{BI86} have proved that the critical points of both functionals
are the Riemannian metrics $h\in {\mathcal A}$  whose Ricci tensor
$\rho$ is Hermitian with respect to the corresponding almost complex
structure $J$, i.e.,
\begin{equation}\label{H-R}
\rho(JX,JY)=\rho(X,Y),\quad X,Y\in TN.
\end{equation}
The K\"ahler metrics satisfy this condition, and Blair and Ianu\c{s}
raised the question of whether a compact almost K\"ahler manifold
with Hermitian Ricci tensor is K\"ahlerian. This question motivated
the following result.

\begin{thm}\label{H-R on Z} \rm (\cite{DM90})
Let $(M,g)$ be a connected oriented real-analytic Riemannian
manifold. If the Ricci tensor of the twistor space $({\mathcal
Z},h_t)$ is
${\mathcal J}_n$-Hermitian, $n=1$ or $n=2$, then either\\
$(i)$ $(M,g)$ is Einstein and self-dual\\
or\\
$(ii)$ $(M,g)$ is self-dual with constant scalar curvature
$s=12/t$ and, for each point of $M$, at least three eigenvalues of
its Ricci operator coincide.

\smallskip

Conversely, if $(M,g)$ is a smooth oriented Riemannian
four-manifold satisfying $(i)$ or $(ii)$, then the Ricci tensor of
$({\mathcal Z},h_t)$ is ${\mathcal J}_n$-Hermitian.
\end{thm}

\smallskip

\noindent {\bf Examples}. \rm (\cite{DM90}) {\bf 1}. Let $M$ be an
Einstein self-dual manifold with negative scalar curvature. Then,
by \cite{M87}, $({\mathcal J}_2,h_t)$ for $t=-12/s$ is an almost
K\"ahler structure on the twistor space ${\mathcal Z}$. This
structure is not K\"ahlerian since, by the Eells-Salamon result
mentioned above, the almost complex structure ${\mathcal J}_2$ is
not integrable. On the other hand, the Ricci tensor of the metric
$h_t$ is ${\mathcal J}_2$-Hermitian by Theorem~\ref{H-R on Z}.
Thus, if $M$ is compact, the twistor space $({\mathcal
Z},h_t,{\mathcal J}_2)$ gives a negative answer to the
Blair-Ianu\c{s} question. Examples of compact Einstein self-dual
manifolds with negative scalar curvature can be found in
\cite{V86}. Multiplying  the twistor space of such a manifold by
K\"ahler manifolds, one can construct examples of non-K\"ahler
almost K\"ahler manifolds of arbitrary even dimension $\geq 6$
which have Hermitian Ricci tensors. Other examples of such
manifolds can be obtained as twistor spaces of quaternionic
K\"ahler manifolds \cite {AGI}. In dimension four the
Blair-Ianu\c{s} problem is not completely solved yet. In this case
positive results under certain additional conditions have been
proved in \cite{AD00,Dra94,Dra95,Dra99} (see also \cite[Sec.
10.2]{Blair})

\smallskip
\noindent {\bf 2}. The Riemannian product $M = S^1\times S^3$ is a
non-Einstein manifold satisfying the conditions (ii) of
Theorem~\ref{H-R on Z}. Other examples of such manifolds can be
obtained as warped-products of $S^1$ and $S^3$, see, for example,
\cite{Der80}.

\smallskip

The twistor space construction can be used to obtain other examples
of almost Hermitian manifolds with Hermitian Ricci tensor. Let
$(M,g,J)$ be a $4$-dimensional almost Hermitian manifold with the
orientation induced by the almost complex structure $J$. Then $J$ is
a section of the positive twistor bundle $\pi:{\mathcal Z}_{+}\to
M$. Taking the horizontal lift of $J$ and the complex structure of
the fibre of ${\mathcal Z}_{+}$ we define an almost complex
structure ${\mathcal J}$ compatible with the metrics $h_t, t>0$.
More precisely, for $\sigma\in{\mathcal Z}$, $X\in
T_{\pi(\sigma)}M$, and $V\in{\mathcal V}_{\sigma}$ we set
$$
{\mathcal J}X^h_{\sigma}=(JX)^h_{\sigma},\quad {\mathcal
J}V=\sigma\times V.
$$
The geometric conditions for integrability of ${\mathcal J}$ have
been obtained in \cite{De} and the Gray-Hervella classes of the
almost Hermitian structure  $(h_t,{\mathcal J})$ have been
determined in \cite{ADM14}.

\begin{thm}\label{Hung Ricci} \rm (\cite{DM18})
The Ricci tensor of the almost Hermitian manifold
$({\mathcal Z}_{+},h_t,{\cal J})$ is Hermitian if and only if the
base manifold $(M,g)$ is Einstein and anti-self-dual.
\end{thm}

\noindent  {\bf Examples}. \rm (\cite{DM18}) According to
\cite[Theorem 1]{ADM14}, the almost Hermitian structure $(h_t,{\cal
J})$ is K\"ahler exactly when $(M,g,J)$ is K\"ahler and Ricci flat.
Thus, in order to construct compact non-K\"ahler twistor spaces
$({\cal Z},h_t,{\cal J})$ with Hermitian Ricci tensor we need
examples of compact, Einstein, anti-self-dual, non-K\"ahler almost
Hermitian manifolds $(M,g,J)$. We consider three cases according to
the sign of the scalar curvature $s$ of such a manifold.

\smallskip

\noindent {\bf 1}. Case $s>0$. In this case, by the Hitchin and
Friedrich-Kurke result we have mentioned, $(M,g)$ is isometric
either to the $4$-sphere $S^4$ with the round metric or to the
complex projective space $\overline{{\mathbb {CP}}^2}$ with the
opposite orientation and the Fubini-Study metric. As is well-known,
none of these manifolds admits an almost complex structure for
topological reasons.

\smallskip

\noindent {\bf 2}. Case $s<0$. As  C. LeBrun pointed out to us,
Conder and Maclachlan \cite{CM05} have constructed a compact
orientable Riemannian manifold $(M,g)$ of constant negative
sectional curvature with Euler characteristic $\chi=16$. The
signature of $M$ is zero by the well-known integral formula
$$
\tau=\frac{1}{12\pi^2}\int\limits_{M}(||{\mathcal
W}_{+}||^2-||{\mathcal W}_{-}||^2)vol_g
$$
since both  half-Weyl tensors ${\cal W}_{\pm}$ vanish. In
particular, the intersection form of $M$ is indefinite. We also have
$\tau+\chi\equiv 0~ mod\, 4$. Hence, by a version of Ehresmann-Wu
theorem due to O. Saeki (see, for example, \cite[Theorem 8
(A)]{Ma91}, $M$ admits an almost complex structure $I$. Then, as is
well-known, $M$ admits also an almost complex structure $J$
compatible with the metric of $M$. Indeed, let $g_0$ be a Riemannian
metric on $M$ compatible with $I$, for example take $g_0(X,Y)=
g(X,Y) + g(IX,IY)$. Define a symmetric positive endomorphism $G$ of
$TM$ by $g_0(GX,Y)=g(X,Y)$. Then $J= G^{-1/2} I G^{1/2}$ is an
almost complex structure compatible with the metric $g$. This almost
complex structure is not integrable as the following lemma shows.
\begin{lemma}\label{ASD NP} \rm (\cite{DM18})
Every compact anti-self-dual Hermitian surface $(M,g,J)$ of
non-positive scalar curvature is K\"ahler and scalar flat.
\end{lemma}
Note also that if a Riemannian manifold $(M,g)$ admits a compatible
almost complex structure $J$, it possesses many such structures
inducing the same orientation as $J$. This can be seen, for example,
by means of the exponential map of the twistor space $({\cal
Z}_{+},h_t)$ of $(M,g)$ endowed with the orientation induced by $J$
\cite{D17, DM18}.

\smallskip

\noindent {\bf 3}. Case $s=0$. It is a result of Hitchin
\cite{Hit74} that every compact, Ricci flat, anti-self-dual,
$4$-manifold is either flat, or is a $K3$-surface, or an Enriques
surface, or the quotient of an Enriques surface by a free
anti-holomorphic involution.  Note also that a K\"ahler surface is
anti-self-dual if and only if it is scalar flat. This well-known
fact follows from the K\"ahler curvature identities which imply that
the eigenvalues of the half-Weyl operator ${\cal W}_{+}$ of a
K\"ahler surface are $s/6,-s/12,-s/12$ (see, for example,
\cite{ApDav00}) .

Now, let $(M,g,J)$ be a compact Ricci flat K\"ahler surface (a
Calabi-Yau surface). Let $J_u(p)=exp_{J(p)}[uV(p)]$ be a
$1$-parameter deformation of the K\"ahler structure $J$,  where $V$
is a non-zero compactly supported section of the pull-back bundle
$J^{\ast}{\mathcal V}\to M$. Then every $J_u$, $u\neq 0$, is
non-integrable. Otherwise, by Lemma~\ref{ASD NP}, $J_u$ would be
K\"ahler and we would have $J_u=J$ everywhere since $J_u=J$ outside
of $supp\, V$. Hence $V=0$ on $M$, a contradiction. Thus $(M,g,J_u)$
is a compact, Ricci flat, anti-self-dual, strictly almost Hermitian
manifold.

Finally, let us note that the twistor spaces $({\mathcal
Z}_{+},h_t,{\mathcal J}_u)$ of the almost Hermitian $4$-manifolds
$(M,g,J_u)$  belong to the Gray-Hervella class ${\mathcal G}_2$
\cite{ADM14}. Moreover, by Proposition~\ref{Hung Ricci}  and
\cite[Theorem 1]{ADM14}, it follows that ${\mathcal G}_2$ is the
only possible Gray-Hervella class of the twistor spaces
$({\mathcal Z}_{+},h_t,{\mathcal J})$ with Hermitian Ricci tensor.

\section{K\"ahler curvature identities on twistor
spaces}\label{KCI} In order to generalize results in K\"ahler
geometry, A. Gray \cite{G76} has introduced three classes of almost
Hermitian manifolds whose curvature tensor resembles that of a
K\"ahler manifold. On an almost Hermitian manifold $(N,h,J)$, these
classes are defined by the following curvature identities:
$$
\begin{array}{l}
{\mathcal A}{\mathcal H}_1:~ R(X,Y,Z,W)=R(X,Y,JZ,JW)\\[8pt]
{\mathcal A}{\mathcal H}_2:~
R(X,Y,Z,W)=R(JX,JY,Z,W)+R(JX,Y,JZ,W)+R(JX,Y,Z,JW)\\[8pt]
{\mathcal A}{\mathcal H}_3:~R(X,Y,Z,W)=R(JX,JY,JZ,JW),
\end{array}
$$
where, as usual, $R(X,Y,Z,W)=h(R(X,Y)Z,W)$ for $X,Y,Z,W\in TN$.
These identities have been used in \cite{TV81} for finding
irreducible components of the space of curvature tensors on an
Hermitian vector space under the action of the unitary group. They
have also been a useful tool for characterizing the K\"ahler
manifolds in various classes of almost Hermitian manifolds, to quote
just a few papers \cite{ AAD02-1, AAD02-2, AD00, S89, Va80, Va82}.
Note that in the last years there has been an intensive study of
Hermitian metrics which are K\"ahler-like in the sense that the
curvature tensor of either the Levi-Civita, Chern,
Bismut-Str\"ominger, or more generally, a Gauduchon connection
\cite{Gaud}  has the same symmetries as the curvature tensor of a
K\"ahler metric(see, for example, \cite{AUV,YZ,YZZ}).

The next theorem gives geometric characterizations of the oriented
Riemannian four-manifolds $(M,g)$ whose twistor spaces $({\mathcal
Z},h_t,{\mathcal J}_n), n=1,2$, belong to one of the K\"ahler
curvature classes listed above.

\begin{thm}\label{KCI} \rm (\cite{DMG95})
$(i)$ $({\mathcal Z},h_t,{\mathcal J}_n)\in {\mathcal A}{\mathcal
H}_3$ is equivalent to $({\mathcal Z},h_t,{\mathcal J}_n)\in
{\mathcal A}{\mathcal H}_2$ and holds if and only if $(M,g)$ is
Einstein and self-dual, $n=1$ or $2$.

$(ii)$ $({\mathcal Z},h_t,{\mathcal J}_1)\in {\mathcal A}{\mathcal
H}_1$ if and only if $(M,g)$ is Einstein and self-dual with scalar
curvature $s=0$ or $s=12/t$.

$(iii)$ $({\mathcal Z},h_t,{\mathcal J}_2)\in {\mathcal A}{\mathcal
H}_1$ if and only if $(M,g)$ is Einstein and self-dual with $s=0$,
\end{thm}

\smallskip

\noindent {\bf Remarks}. \rm (\cite{DMG95})  {\bf 1}. By a result of
S. Goldberg \cite{Gol}, every compact almost K\"ahler manifold of
class ${\mathcal A}{\mathcal H}_1$ is K\"ahler, and A. Gray
\cite[Theorem 5.3]{G76} has raised the question of whether the same
is true under the weaker condition ${\mathcal A}{\mathcal H}_2$.
Now, let $(M,g)$ be an Einstein self-dual four-manifold with
negative scalar curvature $s$. For $t=-12/s$, the twistor space
$({\mathcal Z},h_t,{\mathcal J}_2)$ is an almost K\"ahler manifold
of class ${\mathcal A}{\mathcal H}_2$ by \cite{M87} and
Theorem~\ref{KCI}. This manifold is not K\"ahler, since the almost
complex structure ${\mathcal J}_2$ is never integrable. So, we have
a negative answer to  Gray's question.

\smallskip

\noindent {\bf 2}. Let $(M,g)$ be a Ricci-flat self-dual
four-manifold. Then $({\mathcal Z},h_t,{\mathcal J}_2)$, $t>0$, is a
quasi K\"ahler manifold \cite{M87} of class ${\mathcal A}{\mathcal
H}_1$ which is not K\"ahler. Thus, the Goldberg result cannot be
extended to quasi K\"ahler manifolds. In the case when $M={\mathbb
R}^4$ the twistor space is ${\mathcal Z}={\mathbb R}^4\times S^2$
and we recover an example of A. Gray \cite{G76}.

By a result of I. Vaisman \cite{Va82}, every compact Hermitian
surface of class ${\mathcal A}{\mathcal H}_1$  is K\"ahler. The
twistor space $({\mathcal Z},h_t,{\mathcal J}_1)$ is a
non-K\"ahler Hermitian manifold of complex dimension $3$ and of
class ${\mathcal A}{\mathcal H}_1$ by \cite{FK82},\cite{AHS} and
Theorem~\ref{KCI}. If $M$ is compact, then ${\mathcal Z}$ is also
compact, and we see that the Vaisman result is not true in complex
dimensions greater than $2$.

\smallskip

\noindent {\bf 3}. If $(M,g)$ is an Einstein  self-dual
four-manifold with scalar curvature $s>0$, then ${(\mathcal
Z},h_t,{\mathcal J}_1), t=12/s,$ is a K\"ahler manifold \cite{FK82}
and hence of class ${\mathcal A}{\mathcal H}_1$. In fact, in this
case, as we have already mentioned, either $M=S^4$ or $M={\mathbb
C}{\mathbb P}^2$, so $({\mathcal Z},h_t, {\mathcal J}_1)$ is either
${\mathbb C}{\mathbb P}^3$ or the complex flag manifold $F_{1,2}$
with their standard K\"ahler structures.

\section{$\ast$-Einstein twistor spaces} It is well known that the
Ricci/Chern form of a K\"ahler manifold is the image
$\mathcal{R}(\Omega)$ of the K\"ahler form $\Omega$ under the action
of the curvature operator $\mathcal{R} \in End(\Lambda^2)$. For an
arbitrary almost Hermitian manifold $(N,h,J )$, the $2$-form
$\mathcal{R}(\Omega)$ is neither closed nor of type $(1,1)$, but it
is still closely related to the Ricci form of the canonical
Hermitian connection which represents the first Chern class of
$(N,J)$. The tensor $\rho^{\ast}$ associated to
$\mathcal{R}(\Omega)$ by $$\rho^{\ast}(X,Y) =
\mathcal{R}(\Omega)(X,J Y)= Trace(Z\to R(JZ,X)JY)$$ has been
introduced by S. Tachibana \cite{Ta59}, and is known in the
literature as the $\ast$-Ricci tensor. This tensor then appeared in
almost Hermitian geometry in different contexts. For example, it has
been used by A. Gray \cite{Gr76-2} for studying nearly K\"ahler
manifolds and by F. Tricceri and L. Vanhecke \cite{TV81} for
describing the irreducible components of the space of curvature
tensors on a Hermitian vector space under the action of the unitary
group. The $\ast$-Ricci tensor also plays an important role in the
theory of harmonic almost complex structures, developed recently by
C. Wood \cite{Wo95}.

An almost Hermitian manifold is said to be weakly $\ast$-Einstein if
its $\ast$-Ricci tensor is a multiple of the metric, i.e. if the
K\"ahler form is an eigenvector of the curvature operator. Unlike
K\"ahler-Einstein manifolds, the multiple (usually called the
$\ast$-scalar curvature) need not be a constant and when this holds
the manifold is called $\ast$-Einstein. As we have already
mentioned, for K\"ahler manifolds the Einstein and weakly
$\ast$-Einstein conditions coincide, so it is natural to ask whether
there is a relation between them for more general almost Hermitian
manifolds. The curvature decomposition (\ref{dec}) implies that in
real dimension four the weakly $\ast$-Einstein condition holds if
and only if the traceless Ricci tensor is $J$-anti-invariant and the
K\"ahler form is an eigenvector of the self-dual Weyl operator
$\mathcal{W}_+$. Since, for a Hermitian $4$-manifold, the latter
condition is equivalent to $\mathcal{W}_+$ being degenerate (see
\cite{AG}), it follows from the Riemannian Goldberg-Sachs theorem
\cite{AG,Nu,PB} that any Einstein Hermitian metric is weakly
$\ast$-Einstein. For almost K\"ahler 4-manifolds, it is still an
open question whether the Einstein condition implies the weakly
$\ast$-Einstein one, although J. Armstrong \cite{Arm} has explicitly
described all weakly $\ast$-Einstein strictly almost K\"ahler
Einstein 4-manifolds. This, combined with a result of K. Sekigawa
\cite{Sek}, shows that such manifolds can never be compact, so the
positive answer to the question above would imply the well-known
Goldberg conjecture \cite{Gol} that any compact almost K\"ahler
Einstein 4-manifold must be K\"ahler. In higher dimensions, the
(weakly) $\ast$-Einstein condition has not been so well studied and
it seems that the main reason for that is the lack of interesting
examples. Because of that in \cite {DGM01} the authors studied the
twistor spaces of oriented Riemannian 4-manifolds as a source of
$6$-dimensional examples of $\ast$-Einstein almost Hermitian
manifolds and showed that some four-dimensional results on the
$\ast$-Einstein condition cannot be extended to higher dimensions.

The $\ast$-Ricci tensor $\rho^{\ast}_{t,n}$ of the twistor space
$({\mathcal Z},h_{t},{\mathcal J}_n)$, $t>0, n=1,2$, can be computed
in terms of the curvature of the base manifold $(M,g)$ using the
formula for the sectional curvature of $({\mathcal Z},h_{t})$ in
Proposition~\ref{Sec} and the well-known expression of the
Riemannian curvature tensor by means of sectional curvatures.

\begin{prop}\label{Prop6} \rm (\cite{DGM01}) Let $E,F\in T_{\sigma }{\cal Z}$
 and $X=\pi _{*}E$, $Y=\pi _{*}F$, $A={\cal V}E$, $B={\cal V}F$. Then
$$
\begin{array}{ll}
\rho^{\ast}_{t,n}(E,F)&= [1 +(-1)^{n+1}]g({\cal R}(\sigma ),%
X\land K_{\sigma }Y)-\displaystyle{\frac{t}{2}}g(R(X\land K_{\sigma
}Y)\sigma,R(\sigma)\sigma)
\\[6pt]
&+ \displaystyle{\frac{t}{4}}\mbox{ Trace}(Z\to g(R(X\land
Z)\sigma,R(K_{\sigma }Z\land
K_{\sigma}Y)\sigma))\\[6pt]
&+ (-1)^{n+1}\mbox{ Trace}({\cal V}_{\sigma}\ni C\to
g(R(C)X,R(\sigma\times
C)K_{\sigma }Y))\\[6pt]
&+ \displaystyle{\frac{t}{2}} (-1)^{n}g((\nabla _{X}{\cal
R})(\sigma),B)+\displaystyle{\frac{t}{2}}
g((\nabla _{K_{\sigma }Y}{\cal R})(\sigma ),\sigma\times A)\\[8pt]
&+ [1
+(-1)^{n+1}tg({\cal R}(\sigma ),\sigma
)]g(A,B)\\[6pt]
&+ (-1)^{n+1}\displaystyle{\frac{t^2}{4}}\mbox{ Trace}(Z\to
g(R(\sigma\times A)K_{\sigma }Z,R(B)Z)),
\end{array}
$$
where $K_{\sigma }$ is the complex structure on $T_{\pi(\sigma)}M$
determined  by $\sigma$.
\end{prop}

In the case when the base manifold $(M,g)$ is Einstein and
self-dual the formula for $\rho^{\ast}_{t,n}$ simplifies
significantly:

\begin{cor}\label{Cor1} Let $(M,g)$ be an Einstein self-dual $4$-manifold with
scalar curvature $s$. Then

$$
\rho^{\ast}_{t,n}(E,F) = \frac{1}{12}[ (1 + (-1)^{n+1})s +
\frac{t}{24} (1 + (-1)^n)s^2]g(X,Y)+$$
 $$+[1 + (-1)^{n+1}\frac{ts}{6}+ (-1)^n (\frac{ts}{12})^2]g(A,B).$$
\end{cor}
The above formulas can be used to obtain the following geometric
characterization of the $\ast $-Einstein twistor spaces.

\begin{thm}\label{Star-Ein on Z}\rm (\cite{DGM01})
Let $(M,g)$ be an oriented Riemannian $4$-manifold with scalar
curvature $s$.

\noindent $(i)$ The twistor space $({\mathcal Z},h_t,{\mathcal
J}_1)$ is $\ast$-Einstein if and only if $(M,g)$ is Einstein,
self-dual and $t|s|=12$,

\noindent $(ii)$ The twistor space $({\mathcal Z},h_t,{\mathcal
J}_2)$ is $\ast$-Einstein if and only if $(M,g)$ is Einstein,
self-dual and $ts=6$ .
\end{thm}

A crucial role in the proof of Theorem \ref{Star-Ein on Z} is played
by the following result essentially due to C. LeBrun and V.
Apostolov (private communications, 2000) which is also of
independent interest.

\begin{lemma}\label{lemmaLA} (\cite{DGM01}) There is no self-dual
manifold $(M,g)$ whose Ricci operator has constant eigenvalues
$(\lambda,\mu,\mu,\mu)$ with $\lambda\neq 0$ and $\lambda\neq\mu$.
\end{lemma}
\smallskip

\noindent {\bf Remarks}. \rm (\cite{DGM01})  \noindent {\bf 1}. A
Hermitian metric on a compact complex surface $(M,J)$ is
$\ast$-Einstein if and only if it is locally conformally K\"ahler
and the traceless Ricci tensor is $J$-anti-invariant \cite{AG 97}.
In higher dimensions however the $\ast$-Einstein condition does not
imply any of these two properties, as can be seen by considering the
twistor space $({\cal Z},h_t,{\mathcal J}_1)$  of a compact
self-dual Einstein manifold $(M,g)$ with negative scalar curvature
$s$ and $t= -12/s$. By Theorem~\ref{Star-Ein on Z}, the
6-dimensional Hermitian manifold $({\cal Z},h_t,{\mathcal J}_1)$ is
$\ast$-Einstein, but is neither locally conformally K\"ahler
\cite{M87}, nor with ${\mathcal J}_1$-anti-invariant traceless Ricci
tensor \cite{DM90}.

\smallskip

\noindent {\bf 2}. By a result of V. Apostolov \cite{Ap95}, any
compact $\ast$-Einstein Hermitian surface of negative $\ast$-scalar
curvature is K\"ahler. The twistorial example above shows that the
analogous statement is false in higher dimensions.

\smallskip

\noindent {\bf 3}. Recall that the twistor space $({\mathcal
Z},h_t)$ is an Einstein manifold if and only  the base manifold
$M$ is Einstein  and self-dual  with positive scalar curvature
$s=6/t$ or $ s=12/t$. Thus $({\mathcal Z},h_t,{\mathcal J}_1),
t=s/6$,  is an Einstein Hermitian manifold of real dimension $6$
which is neither locally conformally K\"ahler \cite{M87} nor
$\ast$-Einstein (Theorem~\ref{Star-Ein on Z}). Recall also that if
$M=S^{4}$ or $M={\Bbb C}{\Bbb P}^2$, then ${\mathcal Z}= {\Bbb
C}{\Bbb P}^3$ or ${\mathcal Z}=F_{1,2}=SU(3)/S(U(1)\times
U(1)\times U(1))$, and $(h_t,{\mathcal J}_1)$ for $t=12/s$ is  the
standard K\"ahler-Einstein structure on ${\mathcal Z}$. For
$t=6/s$, $({\mathcal Z},h_t)$ is a Riemannian $3$-symmetric space
\cite{WG68} and ${\mathcal J}_2$ is its canonical almost complex
structure. In this case $({\mathcal Z},h_t,{\mathcal J}_2)$ is a
$\ast$-Einstein nearly K\"ahler manifold by a result of Gray
\cite{Gr76-2}. Note also that for $M=S^4$ and $t=6/s$, $h_t$ is
the "squashed" Einstein metric on ${\Bbb C}{\Bbb P}^3$
\cite[Example 9.83]{Besse}.

\section{Curvature properties of the Chern connection on twistor
spaces}

 It is well-known \cite{Lich, Gaud} that every almost
Hermitian manifold admits a unique connection for which the almost
complex structure and the metric are parallel, and the $(1,1)$-part
of the torsion vanishes. It is usually called the Chern connection
because, in the integrable case, it coincides with the Chern
connection \cite{Chern} of the tangent bundle considered as a
Hermitian holomorphic bundle. This connection plays an important
role in (almost) complex geometry since, by the Chern-Weil theory,
the Chern classes of the manifold are directly related to its
curvature. Note also the classification result of Boothby
\cite{Boot} who proved  that the compact Hermitian manifolds with
flat Chern connection are exactly the quotients of complex Lie
groups equipped with left invariant Hermitian metrics.

Motivated by the works of S. Donaldson \cite{Don} and C. LeBrun
\cite{LB}, V. Apostolov and T. Dragichi \cite{AD} have proposed to
study the problem of existence of almost K\"ahler structures of
constant Hermitian scalar curvature and/or type $(1,1)$ Ricci form
of its Chern connection (from now on we refer to it as the first
Chern form). One of our goals in \cite{DGM02} was to show that the
twistor space of any self-dual Einstein $4$-manifold of negative
scalar curvature admits such an almost K\"ahler structure.

Given an almost Hermitian manifold $(N,g,J)$, denote by $\nabla $
the Levi-Civita connection of $h$. Then the Chern connection
$\nabla^c$  of $(N,g,J)$ is defined by (see, for example,
\cite[Theorem 6.1]{GBN}):
\begin{equation}\label{eq 5.1}
\begin{array}{c}
g(\nabla^c_XY,Z)=g(\nabla_XY,Z)+\displaystyle{\frac{1}{2}}g((\nabla_XJ)(JY),Z)\\[6pt]
+\displaystyle{\frac{1}{4}}g((\nabla_ZJ)(JY)-(\nabla_YJ)(JZ)-(\nabla_{JZ}J)(Y)+(\nabla_{JY}J)(Z),X)
\end{array}
\end{equation}
It belongs to the distinguished 1-parameter family of Hermitian
connections $\nabla^u, u\in \mathbb{R}$, defined by P. Gauduchon
\cite{Gaud} :

\begin{equation}
\begin{array}{c}
g(\nabla^u_XY,Z)=g(\nabla_XY,Z)+\displaystyle{\frac{1}{2}}g((\nabla_XJ)(JY),Z)\\[6pt]
+\displaystyle{\frac{u}{4}}g((\nabla_ZJ)(JY)-(\nabla_YJ)(JZ)-(\nabla_{JZ}J)(Y)+(\nabla_{JY}J)(Z),X)
\end{array}
\end{equation}
The Chern connection corresponds to $u=1$, whereas for $u=-1$ we
obtain the so-called Bismut (or Str\"ominger) connection
\cite{Str,Bis}.

Let $\Omega$ be the K\"ahler $2$-form of $(N,g,J)$ and $\delta\Omega$
the co-differential of $\Omega$ with respect to $\nabla$. Denote
by $\varphi$ and $\psi$ the 2-forms on $N$ defined by
\begin{equation}\label{eq 5.2}
\varphi(X,Y)=Trace(Z\to g((\nabla_XJ)(JZ),(\nabla _YJ)(Z)))
\end{equation}
\begin{equation}\label{eq 5.3}
\psi(X,Y)=\rho^*(X,JY)
\end{equation}
where $\rho^*$ is the $*$-Ricci tensor of $(N,g,J)$.

The formula in the next lemma appears in \cite{Gaud} without proof
and we refer the reader to \cite{DGM02} for its proof.

\begin{lemma}\label{Lemma 5.1} The first Chern form $\gamma^u$ of the Gauduchon connection $\nabla^u$
on an almost Hermitian manifold $(N,g,J)$ is given by
$$
 8\pi\gamma^u=-\varphi-4\psi+2ud\delta\Omega
$$
\end{lemma}

Let $(M,g)$ be an oriented Riemannian 4-manifold with twistor space
${\mathcal Z}$. Denote by $D^c_{n}$ the Chern connection of the
almost-Hermitian manifold $({\mathcal Z}$, $h_{t}$,
$\mathcal{J}_{n})$, $n=1,2$, and by $\gamma_{t,n}$ its first Chern
form. In the case when the base manifold $(M,g)$ is self-dual, an
explicit formula for $\gamma_{t,1}$ has been given by P. Gauduchon
\cite{G}. For an arbitrary oriented Riemannian 4-manifold $(M,g)$,
the first Chern forms $\gamma_{t,n}$, $n=1,2$,  can be computed by
means of the following formula.

\begin{prop}\label{Proposition 5.3} {\rm (\cite{DGM02})}The first Chern form $\gamma_{t,n}$ of the twistor space $({\mathcal
Z},h_t,\cal{J}_n)$, $n=1,2,$ is given by
$$
2\pi\gamma_{t,n}(E,F)=[1+(-1)^{n+1}][g({\cal R}(\sigma),X\land Y)+%
g(A,\sigma\times B)]
$$
where $E,F\in T_{\sigma}{\cal Z}$  and $X=\pi_{*}E$, $Y=\pi_{*}F$,
$A={\cal V}E$, $B={\cal V}F$.
\end{prop}

Now,  we consider the problem of when the curvature tensor $R^c_n$
of the Chern connection $D^c_n$ is of type $(1,1)$, i.e.
$$R^c_n(J_nE,J_nF)G=R^c_n(E,F)G$$ for all $E,F,G\in T{\cal Z}$.

\begin{prop}\label{Proposition 5.4} \rm (\cite{DGM02}) $(i)$ The curvature tensor $R^c_1$ is of type $(1,1)$ if and
only if the base manifold $(M,g)$ is self-dual.

$(ii)$ The curvature tensor $R^c_2$ is of type $(1,1)$ if and only
if the base manifold $(M,g)$ is Einstein and self-dual.

\end{prop}

The next proposition solves the problem when of the Chern
connections $D^c_{1}$ and $D^c_{2}$ of a twistor space have constant
holomorphic sectional curvatures.

\begin{prop}\label{Proposition 5.5} \rm (\cite{DGM02}) $(i)$ The Chern connection $D^c_{1}$ of the almost-Hermitian manifold
$({\cal Z},h_{t},\cal{J}_{1})$ has constant holomorphic sectional
curvature $\kappa$ if and only if $\kappa>0$, the base manifold
$(M,g)$ is of constant sectional curvature $\kappa$, and
$t=1/\kappa$.

$(ii)$  The holomorphic sectional curvature of the Chern
connection $D^c_{2}$ of $({\cal Z},h_{t},\cal{J}_{2})$ is never
constant.
\end{prop}

\section{Holomorphic curvatures of twistor spaces}

Given an almost Hermitian manifold $(M,g,J)$ one can define various
types of curvatures related to the almost Hermitian structure
$(g,J)$. The most important are the holomorphic sectional curvature
\cite {KN} and the holomorphic, Hermitian, and orthogonal (totally
real) bisectional curvatures \cite {GK}, \cite {Ba}, \cite {BIG}.
These curvatures have intensively been studied on K\"ahler manifolds
and a lot of important results have been obtained. For example, the
well-known uniformization theorem for complete K\"ahler manifolds of
constant holomorphic sectional curvature states that any such
manifold is either a complex projective space ${\Bbb C}{\Bbb P}^n$
with the Fubini-Study metric, a quotient of ${\Bbb C}^n$ with the
flat metric or a quotient of the unit ball in ${\Bbb C}^n$ with the
hyperbolic metric \cite {KN}. Moreover, by the solution of the
Frankel conjecture given by Mori \cite {Mor} and by Siu and Yau
\cite {SY}, we know that the complex projective spaces are the only
compact complex manifolds admitting K\"ahler metrics of positive
holomorphic bisectional curvature. Note also that Mok \cite {Mok}
has proved the so-called generalized Frankel conjecture stating that
any compact simply-connected K\"ahler manifold with nonnegative
holomorphic bisectional curvature is biholomorphic to a compact
Hermitian symmetric space. We refer the reader to \cite{KS}, \cite
{CX}, \cite{GZ} for analogous results under some weaker conditions
on the holomorphic bisectional curvature. The case of negative
holomorphic bisectional curvature is not so rigid. For example,
recently To and Yeung \cite{TY} have constructed such K\"ahler
metrics on any Kodaira surface.

In the non-K\"ahler case the holomorphic curvatures mentioned
above are not so well studied. Complete results have been obtained
only for complex dimension $2$ in which case it has been proved
that every compact Hermitian surface of constant holomorphic or
Hermitian sectional curvature is a complex space form \cite {ADM}.
In higher dimensions it is still an open question posed by Balas
and Gauduchon \cite {Ba,BG}  whether there are compact
non-K\"ahler Hermitian manifolds of non-zero constant holomorphic
sectional curvature of the Chern connection.

\subsection{Holomorphic bisectional curvature}
The holomorphic bisectional curvature $H_{t,n}$ of the twistor
space $({\cal Z},h_{t},\cal{J}_n), n=1,2$, of an oriented
Riemannian $4$-manifold $(M,g)$ can be computed by means of
Proposition~\ref{Sec}. For the sake of simplicity we give the
respective formula only in the case when the base manifold is
self-dual and Einstein.

\begin{prop}\label {Pr1} \rm (\cite{ADM15}) Let $(M,g)$ be a self-dual Einstein manifold with scalar curvature $s$
and let $E,F\in T_{\sigma }{\cal Z}$ be arbitrary $h_t$-unit
tangent vectors with $X=\pi _{*}E$, $Y=\pi _{*}F$, $V={\cal V}E$,
$W={\cal V}F$. Then
\begin{eqnarray}\label{curv}
H_{t,n}(E,F)&=& R(X,K_\sigma X,Y,K_\sigma Y) +t\parallel V\parallel^2\parallel W\parallel^2\nonumber\\
&+&2t(\frac{s}{24})^2(\parallel X\parallel^2\parallel Y\parallel^2-g(X,Y)^2-g(K_\sigma X,Y)^2)\nonumber\\
&+&(-1)^n(2(\frac{ts}{24})^2-\frac{ts}{12})(\parallel
X\parallel^2\parallel W\parallel^2+\parallel Y\parallel^2\parallel
V\parallel^2)\nonumber\\
&+&(2(\frac{ts}{24})^2(1+(-1)^n)-\frac{ts}{12})(g(K_\sigma
X,Y)g(\sigma\times
V,W)\nonumber\\
&+&(-1)^ng(X,Y)g(V,W)),\nonumber\\
\end{eqnarray}
where $K_{\sigma }$ is the complex structure on $T_{\pi(\sigma)}M$
determined  by $\sigma$.
\end{prop}

We next consider two particular cases of Proposition~\ref{Pr1}.

\begin{cor}\label{Cor4.2}
Let $(M,g)$ be a 4-manifold of constant sectional curvature and
scalar curvature $s$. Then
\begin{eqnarray}\label{concurv}
H_{t,n}(E,F)&=&\frac{s}{12}(g(X,Y)^2+g(K_\sigma X,Y)^2)+ t\parallel V\parallel^2\parallel W\parallel^2\nonumber\\
&+& 2t(\frac{s}{24})^2(\parallel X\parallel^2\parallel
Y\parallel^2
-g(X,Y)^2-g(K_\sigma X,Y))\nonumber\\
&+&(-1)^n(2(\frac{ts}{24})^2-\frac{ts}{12})(\parallel
X\parallel^2\parallel W\parallel^2+\parallel Y\parallel^2\parallel
V\parallel^2)\nonumber\\
&+&(2(\frac{ts}{24})^2(1+(-1)^n)-\frac{ts}{12})(g(K_\sigma
X,Y)g(\sigma\times
V,W)\nonumber\\
&+&(-1)^ng(X,Y)g(V,W)).\nonumber\\
\end{eqnarray}

\end{cor}

\smallskip
\begin{cor}\label{Cor4.3} Let $(M,g)$ be a self-dual Einstein manifold with sectional curvature K and scalar curvature
$s$, and let $E\in T_{\sigma }{\cal Z}$ be arbitrary $h_t$-unit
tangent vector with $X=\pi _{*}E$ and $V={\cal V}E$. The
holomorphic sectional curvature of $({\cal Z},h_t,\cal{J}_n)$ is given
by
$$H_{t,n}(E)=K(X,K_\sigma X)\|X\|^{4}+t\|V\|^{4}+(2(\frac{st}{24})^{2}(3(-1)^{n}+1)+(-1)^{n+1}\frac{st}{24})\|X\|^{2}\|V\|^{2}$$\\
\end{cor}

Using Proposition~\ref{Sec} and Corollary~\ref{Cor4.3} we obtain the
following.

\begin{thm}\label{thmHS}\rm (\cite{DM}) $(i)$ The almost Hermitian manifold $({\cal
Z},h_t,{\cal J}_1)$ has constant holomorphic sectional curvature
$\mathcal{X}$ if and only if the base manifold $(M,g)$ has
constant sectional curvature $\mathcal{X}=1/t$.

$(ii)$ The holomorphic sectional curvature of $({\cal Z},h_t,{\cal J}_2)$
is never constant.

\end{thm}

This together with Corollary~\ref{Cor4.2} implies

\begin{thm}\label{thHBS1} \rm (\cite{ADM15})
The holomorphic bisectional curvature of the twistor space $({\cal
Z},h_t,\cal{J}_n), n=1, 2,$ of an oriented Riemannian 4-manifold
$(M,g)$ is never constant.
\end{thm}

In the next theorem, we consider the case when the base manifold
$(M,g)$ is a real space form and  determine all $t>0$ for which the
holomorphic bisectional curvature of its twistor space $({\cal
Z},h_t,J_n)$ is strictly positive. In particular, it follows that
the "squashed" metric on ${\Bbb C\Bbb P}^3$ (\cite{Besse}, Example
9.83) is a non-K\"ahler Hermitian-Einstein metric of positive
holomorphic bisectional curvature. This shows that a recent result
of Kalafat and Koca \cite{KK} in dimension four can not be extended
to higher dimensions.

\begin{thm}\label{th HBS2} \rm (\cite{ADM15}) Let $(M,g)$ be an oriented
Riemannian $4$-manifold of constant sectional curvature.

(i) The holomorphic bisectional curvature of $({\cal Z},h_t,\cal{J}_1)$
is positive if and only if $0<ts<24$.

(ii) If $(M,g)$ is a flat manifold, the holomorphic bisectional
curvature of $({\cal Z},h_t,\cal{J}_n)$ is non-negative, $n=1,2$.
\end{thm}

As an example illustrating Theorem \ref{th HBS2} (ii), let us
consider the twistor space $({\cal Z},h_1,\cal{J}_1)$ of a 4-torus
$T$ with its standard flat metric. Then $Z=T\times S^2$, $h_1$ is
the product metric and $\cal{J}_1$ is the complex structure
defined by Blanchard \cite {Bl}. So, the holomorphic bisectional
curvature of $(T\times S^2,h_1,{\mathcal J_1})$ is non-negative.
Note that $\cal{J}_1$ is not a product of complex structures on
$T$ and $S^2$.

\subsection{Orthogonal bisectional curvature}

The orthogonal (totally real) bisectional curvature $B$ of an almost
Hermitian manifold $(N,h,J)$ is defined in \cite{BIG} by
$$B(X,Y)=h(R(X, JX)Y,JY)$$ for $X,Y\in TN$ such that $X\bot \{Y, JY\}$ and
$||X||=||Y||=1$.  It is well known \cite{H} that the orthogonal
bisectional curvature of  a K\"ahler manifold of complex dimension
$\geq3$ is constant if and only if the holomorphic sectional curvature is constant. So,
it is natural to ask if the same holds for other classes of almost
Hermitian manifolds. The next theorem shows  that this is true for the
twistor spaces of self-dual Einstein 4-manifolds.

\begin{thm}\label{th3} \rm (\cite{ADM15})
Let $(M,g)$ be a self-dual Einstein $4$-manifold. Then its twistor space $(Z, h_{t},
\cal{J}_{n})$ has constant orthogonal bisectional curvature if and only
if $n=1$ and $(M,g)$ is of constant sectional curvature
$\displaystyle \chi=1/t$.
\end{thm}

\noindent {\bf Remark}. Let $(M,g)$ have a constant sectional
curvature. Then the orthogonal bisectional curvature $B_{t,1}$ of
the twistor space $({\cal Z},h_t,\cal{J}_1)$ is strictly positive if
and only if $0<ts<24$.

\subsection{Hermitian bisectional curvature} The Hermitian bisectional
curvature $H^c$ of an almost Hermitian manifold $(N,h,J)$ is defined
as the holomorphic bisectional curvature of its Chern connection. As
we have already noted, the curvature of this connection is directly
related to the Chern classes of $(N,J)$. In particular, if $\gamma$
is the first Chern form of $(N,h,J)$, then for any $X\in TN$ we have

\begin{equation} \label{Chern}
\gamma (X,JX) = \sum_{i=1}^n h(H^c(X,JX)E_i,JE_i),
\end{equation}
where $(E_1,\dots,E_n,JE_1,\dots,JE_n)$ is a unitary frame.

According to Theorem~\ref{thHBS1}, the holomorphic bisectional
curvature of the twistor space of an oriented Riemannian
$4$-manifold is never constant. As for the Hermitian bisectional
curvature, we have the following more general result which was
pointed out to us by S. Kobayashi (private communication, April
2012).

\begin{thm}\label{th5} \rm (\cite{ADM15}) The Hermitian bisectional curvature of a Hermitian
manifold of complex dimension $\geq2$ is never a non-zero
constant.
\end{thm}

\noindent The proof of this theorem uses formula (\ref {Chern})
for the first Chern form $\gamma$ which implies that if the
Hermitian bisectional curvature of a Hermitian manifold is a
non-zero constant $c$, then the manifold is K\"ahler. Hence it is
a complex space form and the well-known formula for its curvature
\cite{KN} implies that $c=0$, a contradiction. Note also that
Theorem~\ref{th5} gives a partial negative answer to the question
of Balas and Gauduchon \cite {Ba,BG} mentioned at the beginning of
this section.

\bigskip

\noindent \textbf{Remark.} (\cite{ADM15}) Formula (\ref {Chern})
for the first Chern form implies that if an almost Hermitian
manifold has non-zero constant Hermitian bisectional curvature,
then it is an almost K\"ahler manifold, i.e. its K\"ahler $2$-form
is closed. Hence it is natural to ask the following questions:

\begin{itemize}

    \item Are there compact non-K\"ahler and non-flat
Hermitian manifolds of complex dimension $\geq3 $ with vanishing
Hermitian bisectional curvature?

     \item Are there compact non-K\"ahler almost K\"ahler manifolds of constant
Hermitian bisectional curvature?
\end{itemize}
\smallskip
By a result of Vezzoni \cite [Theorem 4.8]{Ve}, if $(N,h,J)$ is an
almost K\"ahler manifold whose holomorphic and Hermitian
bisectional curvatures coincide, then it is a K\"ahler manifold.
This result can be extended to a more general class of almost
Hermitian manifolds.

\begin{thm}\label{th41} \rm (\cite{ADM15})
Let $(N,h,J)$ be an almost Hermitian manifold such that
\begin{equation} \label{GrindEQ__4_1_}
({\nabla }_XJ)(X)= \varepsilon ({\nabla }_{JX}J)(JX),
\end{equation}
where $\varepsilon =\pm1$. Then its holomorphic and Hermitian
bisectional curvatures coincide if and only if $(N,h,J)$ is a
K\"ahler manifold.
\end{thm}

\smallskip
\noindent {\bf Remarks}. (\cite{ADM15}) ${\bf 1.}$ According to
the Gray-Hervella terminology \cite{GH} the almost Hermitian
manifolds satisfying (\ref{GrindEQ__4_1_}) with $\varepsilon=1$
are called $G_1$-spaces. This class contains the Hermitian and
nearly K\"ahler manifolds. The identity (\ref{GrindEQ__4_1_}) with
$\varepsilon=-1$ holds for almost K\"ahler  and quasi K\"ahler
manifolds (recall that the quasi K\"ahler  condition is
$(\nabla_XJ)(Y)+(\nabla_{JX}J)(JY)=0$).

\noindent ${\bf 2}.$  The proof of Theorem \ref{th41} shows that
the above mentioned result of Vezzoni for almost K\"ahler
manifolds holds true under the weaker condition that the
holomorphic and Hermitian sectional curvatures coincide.

\smallskip

Finally,  we describe the twistor spaces whose holomorphic and Hermitian
sectional curvatures coincide.

\smallskip

\begin{thm}\label{th51} \rm (\cite{ADM15})
Let $(M,g)$ be an oriented Riemannian $4$-manifold. The holomorphic
and Hermitian sectional curvatures of its twistor space $({\cal
Z},h_{t},\cal{J}_n)$ coincide if and only if $(M,g)$ is a self-dual
Einstein manifold with $ts=12$ for $n=1$ and $ts=6$ for $n=2$.
\end{thm}

\end{document}